\documentclass[12pt,reqno]{amsart}
\setbox0=\hbox{$+$}
\newdimen\plusheight
\plusheight=\ht0
\def\+{\;\lower\plusheight\hbox{$+$}\;}

\setbox0=\hbox{$-$}
\newdimen\minusheight
\minusheight=\ht0
\def\-{\;\lower\minusheight\hbox{$-$}\;}

\setbox0=\hbox{$\cdots$}
\newdimen\cdotsheight
\cdotsheight=\plusheight
\def\cds{\lower\cdotsheight\hbox{$\cdots$}}

\renewcommand{\(}{\left\(}
\renewcommand{\)}{\right\)}

\def\NI{\noindent}

\renewcommand{\pmod}[1]{\,(\textup{mod}\,#1)}
\numberwithin{equation}{section}
\theoremstyle{plain}
\newtheorem{theorem}{Theorem}[section]
\newtheorem{lemma}[theorem]{Lemma}
\newtheorem{corollary}[theorem]{Corollary}

\begin{document}
	\begin{center}{\bf \large  Some new congruences for generalized overcubic partition function }\end{center}\vskip
	5mm
	\begin{center}
		\centerline{\bf Adam Paksok and Nipen Saikia$^{\ast}$}\vskip2mm
		
		{\it Department of Mathematics, Rajiv Gandhi University,\\ Rono Hills, Doimukh, Arunachal Pradesh, India-791112.\\ E. Mail(s): adam.paksok@rgu.ac.in; nipennak@yahoo.com}\\
		$^\ast$\textit{Corresponding author}\end{center}\vskip3mm
	
	\noindent{\bf Abstract:} Amdeberhan et al. (2024) introduced the notion of a generalized overcubic partition function $\overline a_c (n)$ and proved an infinite family of congruences modulo a prime  $p\ge 3$ and some Ramanujan type congruences.  In this paper, we show that $\overline a_{2^\lambda m+t}(n) \equiv \overline a_t (n) \pmod {2^{\lambda+1}}$, where $\lambda \geq1, m\geq0,$ and  $t\geq1$ are integers. We also prove some new congruences modulo $8$ and $16$ for $\overline a_{2m+1}(n), \overline a_{2m+2}(n), \overline a_{8m+3}(n)$, where $m$ is any non-negative integer.
	
	\vskip
	3mm

	\vskip
	3mm
	
	\NI{\bf Mathematics Subject Classifications:} 11P83. 05A17 \vskip
	3mm

	\section{Introduction}
	 For any complex number $a$ and $q$ (with $|q|<1$), define the $q$-product
	\begin{equation}\label{1}
		(a;q)_{\infty}=\prod_{n=0}(1-aq^n).\end{equation}
	 For brevity, we will use the notation  
	$f_h:=(q^h;q^h)_{\infty}$ for any positive integer $h$.\\ 
	A partition of a positive integer $n$ is a finite non-increasing sequence of positive integers $u_1, u_2, \dots, u_k $ satisfying $\sum_{i=1}^{k}u_i=n$. The integers $u_1, u_2, \dots, u_k $  are called  parts or summands of the partition. For example, the partitions of $n=4$ are given by $4, 3+1, 2+2, 2+1+1,  1+1+1+1$. The number of partition of a non-negative integer $n$ is usually  denoted by $p(n)$ and its  generating function is given by 
	\begin{equation}\label{p}\sum_{n\ge0}^{\infty}p(n)q^n=\dfrac{1}{f_1},\qquad p(0)=1.\end{equation} Ramanujan \cite{bc1, sr2, SR} established following beautiful congruences for the  partition function $p(n)$:
	\begin{equation}
		\label{pc}
		p(5n+4) \equiv0\pmod5,\quad
		p(7n+5)\equiv0 \pmod7,\quad p(11n+6)\equiv0\pmod {11}.\end{equation}
	
	An overpartition of a positive integer $n$ is a partition of $n$ wherein first occurrence of a part may be overlined.  The number of overpartitions of $n$ is usually denoted by $\overline p(n)$.  The overpartition function $\overline p(n)$ is first studied by Corteel and Lovejoy \cite{cor} who also  gave the generating function as
	\begin{equation}\label{o}\sum_{n\ge0}^{\infty}\overline p(n)q^n=\dfrac{f_2}{f_1^2}.\end{equation}	
	For instance, $\overline p(3)=8$ with overpartitions of 3 given by  $3, \overline{3}, 2+1, \overline 2+1, 2+\overline 1, \overline 2+\overline 1, 1+1+1, \overline1+1+1$.
	One may see \cite{hs2, hs1, kim1, kmah, xia1over} and references therein for congruence and other arithmetic properties of overpartition function $\overline{p}(n)$.
	
	Chan \cite{chan} studied the cubic partition of a positive integer $n$ in which even parts appear in two colors. If $a(n)$ denotes the number of cubic partitions of $n$, then the generating function of $a(n)$ is given by
	\begin{equation}\label{c}
		\sum_{n=0}^{\infty}a(n)q^n=\dfrac{1}{f_1f_2}.
	\end{equation}
	For arithmetic properties of $a(n)$ see \cite{chan, chan1, chan2, chan3}.
	
	Kim \cite{kim} studied the overcubic partition function $\overline a(n)$  which represents number of all the overlined versions of cubic partition and the generating function for $\overline a(n)$ is given by 
	\begin{equation}\label{oc}
		\sum_{n=0}^{\infty}\overline a(n)q^n=\dfrac{f_4}{f_1^2f_2}.
	\end{equation}
	Recently, Amdeberhan et al. \cite{ass} considered generalized overcubic partition of $n$ denoted by $\overline a_c (n)$, where $c$ is any positive integer. For instance, $\overline a_2 (3)=12$ with relevent partitions given by $ 3, \bar{3}, 2+1, 2+1, \bar{2}+1, \bar{2}+1, 2+\bar{1}, 2+\bar{1}, \bar{2}+\bar{1}, \bar{2}+\bar{1}, 1+1+1, \bar{1}+1+1 $. The generating function of $\overline a_c (n)$ is given by
	\begin{equation}\label{goc}
	\sum_{n=0}^{\infty}\overline a_c (n)=\dfrac{f_4^{c-1}}{f_1^2f_2^{2c-3}}.
	\end{equation} Amdeberhan et al. \cite{ass} proved an infinite family of congruences modulo a prime  $p\ge 3$ and some Ramanujan type congruences. 
     
     In this paper, we prove some new congruences modulo $8$ and $16$ for $\overline a_{2m+1}(n)$, $\overline a_{2m+2}(n)$ and $\overline a_{8m+3}(n)$, where $m$ is any non-negative integer. To prove our congruences, we use some theta-function and $q$-series identities, which are listed in Section $2$.

	\section{Preliminaries}
	Ramanujan's general theta-function $F\left(x,y\right)$ is defined by
	\begin{equation}\label{eq2a}
		F(x, y)=\sum_{\nu=-\infty}^\infty x^{\nu(\nu+1)/2}y^{\nu(\nu-1)/2},~  |xy|<1.
	\end{equation} 
	The three special  cases of $F\left(x,y\right)$ \cite[p. 36, Entry 22 (i), (ii))]{bcb3} are  the theta-functions $\phi(q)$, $\psi(q)$  and  $F(-q)$ given by
	\begin{equation}\label{e0}
		\phi(q):=F(q, q)=\sum_{\nu=0}^{\infty}q^{\nu^2}=(-q;q^2)^{2}_{\infty}(q^2;q^2)_{\infty}=\dfrac{{f
				_2^5}}{{f_1^2}{f_4^2}} ,
	\end{equation}
	\begin{equation}\label{tp1}
		\hspace{-1cm} \psi(q):=F(q,q^3)=\sum_{\nu\geq 0}q^{\nu(\nu+1)/{2}}=\dfrac{(q^2;q^2)_{\infty}}{(q;q^2)_{\infty}}=\dfrac{{f
				_2^2}}{f_1},
	\end{equation}
	and
	\begin{equation}\label{t3}
		\hspace{-.7cm}	F(-q):=F(-q,-q^2)=\sum_{\nu=-\infty}^{\infty}(-1)^\nu q^{\nu(3\nu+1)/{2}}=f_1.\end{equation}
	Employing elementary $q$-operations, it is easily seen that \begin{equation}\label{phim}
		\phi(-q)=\dfrac{(q;q)^{2}_{\infty}}{(q^2;q^2)_{\infty}}=\dfrac{f_{1}^2}{f_2}.
	\end{equation}
	
	\begin{lemma}\label{iii}\cite[Theorem 2.1]{cui} For any  prime $p>2$, we have
		\begin{equation}\label{w4}
			\psi(q)= \sum_{j=0}^{(p-3)/2}{\it q}^{(j^2+j)/2} F\left( {\it q}^{(p^2+(2j+1)p)/2},{\it q}^{(p^2-(2j+1)p)/2}\right)+{\it q}^{(p^2-1)/8}\psi(q^{p^2}).\end{equation} \\Furthermore, $(j^2+j)/2 \not\equiv(p^2-1)/8 \pmod p~ for ~0 \leq j\leq (p-3)/2.$
	\end{lemma}
	
	\begin{lemma}\cite[Theorem 2.2]{cui} For any prime $p\geq 5$, we have
		\begin{equation}\label{u10} 
			f_1=\sum_{\substack{ k={-(p-1)/2} \\ k \ne {(\pm p-1)/6}}}^{(p-1)/2}(-1)^{k}{q}^{(3k^2+k)/2} F\left( {-q}^{(3p^2+(6k+1)p)/2},{-q}^{(3p^2-(6k+1)p)/2}\right)+(-1)^{{(\pm p-1)/6}}{q}^{(p^2-1)/24}f_{p^2},\end{equation}where
		\begin{equation*}
			\dfrac{\pm p-1}{6}
			= \left\{
			\begin{array}{ll} 
				\dfrac{(p-1)}{6},   
				& if~ p \equiv 1\pmod 6 \\
				\dfrac{(-p-1)}{6} ,
				& if~ p \equiv -1\pmod 6. 
			\end{array}
			\right.
		\end{equation*} \\Furthermore, if $\dfrac{-(p-1)}{2}\leq k \leq\dfrac{(p-1)}{2}~and ~k \neq \dfrac{(\pm p-1)}{6}$~ then\\
		$$\dfrac{3k^2+k}{2}\not\equiv \dfrac{p^2-1}{24} \pmod p.$$
	\end{lemma}
	
	\begin{lemma}\cite{ovyt} We have
		\begin{equation}\label{f_1^2}
			\dfrac{1}{{f_1^2}}=\dfrac{{f_8^5}}{{f_2^5}{f_{16}^2}}+2q\dfrac{{f_4^2}{f_{16}^{2}}}{{f_2^5}{f_8}}.\end{equation}
	\end{lemma}
	
	\begin{lemma}\cite{hs1} We have
		\begin{equation}\label{1/phi} 
			\dfrac{1}{\phi(-q)}=\dfrac{1}{\phi(-q^4)^4}\Big(\phi(q^4)^3+2q\phi(q^4)^2\psi(q^8)+4q^2\phi(q^4
			)\psi(q^8)^2+8q^3\psi(q^8)^3\Big).
		\end{equation}
	\end{lemma} 
	
	\begin{lemma}\cite{bcb3}We have
		\begin{equation}\label{f_1^3}
		f_1^3=f_{25}^3\left(D^{-3}(q^5)-3qD^{-2}(q^5)+5q^3-3q^5D^{2}(q^5)-q^6D^{3}(q^5)\right),
		\end{equation}
		where
		$$D(q)=\dfrac{(q^2;q^5)_\infty(q^3;q^5)_\infty}{(q;q^5)_\infty(q^4;q^5)_\infty}.$$
	\end{lemma}

	To end this section, we record the following congruence which can be easily proved using the binomial theorem: 
	For any positive integer $k, m$ and prime $p$, we have
	\begin{equation}\label{yp}
		{f_{pm}}\equiv {f_m}^p \pmod p,\end{equation}
	and 
	\begin{equation}\label{yp2}
		{f_{pm}^{p^{k-1}}}\equiv {f_m}^{p^k} \pmod {p^k}.\end{equation}

	\section{Congrurences for $\overline a_c (n)$}
	\begin{theorem} If $\lambda\geq1$, $m\geq0$, $n\geq0$ and $t\geq1$ are any integers, then we have
		\begin{equation}\label{3.1}
		\overline a_{2^\lambda m+t}(n) \equiv \overline a_t (n) \pmod {2^{\lambda+1}}.
		\end{equation}
	\end{theorem}
	\begin{proof}
		Setting $c=2^\lambda m+t$ in \eqref{goc}, we obtain
		\begin{equation}\label{3.1a}
			\sum_{n=0}^{\infty}\overline a_{2^\lambda m+t}(n)q^n=\dfrac{f_4^{2^\lambda m+t-1}}{f_1^2f_2^{2^{\lambda+1} m+2t-3}}=\dfrac{f_4^{2^\lambda m}}{f_2^{2^{\lambda+1}m}}\cdot \dfrac{f_4^{t-1}}{f_1^2f_2^{2t-3}}.
		\end{equation}
		Employing \eqref{yp2} with ${m=2, p=2, k=\lambda+1}$, we obtain
		\begin{equation}\label{3.1b}
			\sum_{n=0}^{\infty}\overline a_{2^\lambda m+t}(n)q^n\equiv\dfrac{f_4^{t-1}}{f_1^2f_2^{2t-3}}\pmod{2^{\lambda+1}}.
		\end{equation}
		Employing \eqref{goc}, \eqref{3.1} follows from \eqref{3.1b}.
	\end{proof}
	
	\begin{theorem}For any integers $\lambda\geq0$ and $m\geq0$, we have
			\begin{equation}\label{3.2i}
			\overline a_{2m+1}(8\cdot5^{2\alpha}n+5^{2\alpha}) \equiv \overline a_{2m+1} (8n+1) \pmod {8},
		\end{equation}
		\begin{equation}\label{3.2ii}
		\hspace{-1.5cm}	\overline a_{2m+1}(8\cdot5^{2\alpha+1}n+17\cdot5^{2\alpha}) \equiv0 \pmod {8},
		\end{equation}
		and
		\begin{equation}\label{3.2iii}
		\hspace{-1.5cm}	\overline a_{2m+1}(8\cdot5^{2\alpha+1}n+33\cdot5^{2\alpha}) \equiv0 \pmod {8}.
		\end{equation}\end{theorem}
	
	\begin{proof}
		Setting $c=2m+1$ in \eqref{goc}, we obtain
		\begin{equation}\label{3.2a}
			\sum_{n=0}^{\infty}\overline a_{2m+1}(n)q^n=\dfrac{f_4^{2m}f_2}{f_1^2f_2^{4m}}.
		\end{equation}
		Employing \eqref{f_1^2} in \eqref{3.2a}, we obtain
		\begin{equation}\label{3.2b}
				\sum_{n=0}^{\infty}\overline a_{2m+1}(n)q^n=\dfrac{f_4^{2m}f_8^5}{f_2^{4m+4}f_{16}^2}+2q\dfrac{f_4^{2m+2}f_{16}^2}{f_2^{4m+4}f_8}.
		\end{equation}
		Extracting the terms involving $q^{2n+1}$ from \eqref{3.2b}, dividing by $q$ and replacing $q^2$ by $q$, we obtain
		\begin{equation}\label{3.2c}
				\sum_{n=0}^{\infty}\overline a_{2m+1}(2n+1)q^n=2\dfrac{f_2^{2m+2}f_{8}^2}{f_1^{4m+4}f_4}.
		\end{equation}
		Employing \eqref{yp2} with ${k=2, m=1}$ and $p=2$, we obtain
		\begin{equation}\label{3.2d}
				\sum_{n=0}^{\infty}\overline a_{2m+1}(2n+1)q^n\equiv2f_4^3\pmod{8}.
		\end{equation}
		Extracting the terms involving $q^{4n}$ from \eqref{3.2d} and replacing $q^4$ by $q$, we obtain
		\begin{equation}\label{3.2e}
				\sum_{n=0}^{\infty}\overline a_{2m+1}(8n+1)q^n\equiv2f_1^3\pmod{8}.
		\end{equation}
		Employing \eqref{f_1^3} in \eqref{3.2e}, we obtain
		\begin{equation}\label{3.2f}
			\sum_{n=0}^{\infty}\overline a_{2m+1}(8n+1)q^n\equiv2f_{25}^3\left(D^{-3}(q^5)-3qD^{-2}(q^5)+5q^3-3q^5D^{2}(q^5)-q^6D^{3}(q^5)\right)\pmod{8}.
		\end{equation}
			Extracting the terms involving $q^{5n+3}$ from \eqref{3.2f}, dividing by $q^3$ and replacing $q^5$ by $q$, we obtain
		\begin{equation}\label{3.2g}
			\sum_{n=0}^{\infty}\overline a_{2m+1}(8(5n+3)+1)q^n\equiv2f_5^3\pmod{8}.
		\end{equation}
		Extracting the terms involving $q^{5n}$ from \eqref{3.2g} and replacing $q^5$ by $q$, we obtain
		\begin{equation}\label{3.2h}
			\sum_{n=0}^{\infty}\overline a_{2m+1}(8(5^2n+3)+1)q^n\equiv2f_1^3\pmod{8}.
		\end{equation}
		Employing \eqref{3.2e} in \eqref{3.2h}, we obtain
		\begin{equation}\label{3.2j}
			\sum_{n=0}^{\infty}\overline a_{2m+1}(8(5^2n+3)+1)q^n\equiv	\sum_{n=0}^{\infty}\overline a_{2m+1}(8n+1)q^n\pmod{8}.
		\end{equation}
		Iterating \eqref{3.2j} with $n$ replaced by $5^2n+3$ and simplying, we obtain
		\begin{equation}\label{3.3k}
			\sum_{n=0}^{\infty}\overline a_{2m+1}(8\cdot5^{2\alpha}n+5^{2\alpha})q^n\equiv	\sum_{n=0}^{\infty}\overline a_{2m+1}(8n+1)q^n\pmod{8},
		\end{equation}
		for any integer $\alpha\geq0$. Equating the coefficients of $q^n$ on both sides of \eqref{3.3k}, we arrive \eqref{3.2i}.\\
		Combining \eqref{3.2f} and \eqref{3.3k}, we obtain
		\begin{equation}\label{3.2l}
			\sum_{n=0}^{\infty}\overline a_{2m+1}(8\cdot5^{2\alpha}n+5^{2\alpha})q^n\equiv2f_{25}^3\left(D^{-3}(q^5)-3qD^{-2}(q^5)+5q^3-3q^5D^{2}(q^5)-q^6D^{3}(q^5)\right)\pmod{8}.
		\end{equation}
			Extracting the terms involving $q^{5n+2}$ and $q^{5n+4}$ from \eqref{3.2l}, we complete the proof of \eqref{3.2ii} and \eqref{3.2iii}.\end{proof}
	
	\begin{theorem}
		Let  $p\geq5$ be a prime with $\left(\dfrac{-2}{p}\right)=-1$ and $1\leq r \leq p-1$. Then for all integer $m\geq0$ and $\alpha\geq0$, we have
		\begin{equation}\label{3.3}
			\sum_{n=0}^{\infty}\overline a_{2m+2}(8\cdot p^{2\alpha}n+p^{2\alpha})q^n\equiv2f_1f_2\pmod{8},
		\end{equation}and
			\begin{equation}\label{3.3i}
			\overline a_{2m+2}(8\cdot p^{2\alpha+1}(pn+r)+p^{2\alpha+2})\equiv0\pmod{8}.
		\end{equation}\end{theorem}
	
	\begin{proof}
			Setting $c=2$ in \eqref{goc}, we obtain
		\begin{equation}\label{3.3a}
			\sum_{n=0}^{\infty}\overline a_{2}(n)q^n=\dfrac{f_4}{f_1^2f_2}.
		\end{equation}
		Employing \eqref{f_1^2} in \eqref{3.3a}, we obtain
		\begin{equation}\label{3.3b}
				\sum_{n=0}^{\infty}\overline a_{2}(n)q^n=\dfrac{f_4f_8^5}{f_2^6f_{16}^2}+2q\dfrac{f_4^3f_{16}^2}{f_2^{6}f_8}.
		\end{equation}
			Extracting the terms involving $q^{2n+1}$ from \eqref{3.3b}, dividing by $q$ and replacing $q^2$ by $q$, we obtain
		\begin{equation}\label{3.3c}
			\sum_{n=0}^{\infty}\overline a_{2}(2n+1)q^n=2\dfrac{f_2^3f_{8}^2}{f_1^6f_4}.
		\end{equation}
		Employing \eqref{yp2} and then applying \eqref{3.1}, we obtain
		\begin{equation}\label{3.3d}
			\sum_{n=0}^{\infty}\overline a_{2m+2}(2n+1)q^n\equiv2\dfrac{f_2f_4^3}{f_1^2}\pmod{8}.
		\end{equation}
		Again, employing \eqref{f_1^2} in \eqref{3.3d}, we obtain
		\begin{equation}\label{3.3e}
				\sum_{n=0}^{\infty}\overline a_{2m+2}(2n+1)q^n\equiv2\dfrac{f_4^3f_8^5}{f_2^4f_{16}^2}+4q\dfrac{f_4^5f_{16}^2}{f_2^4f_8}\pmod{8}.
		\end{equation}
			Extracting the terms involving $q^{2n}$ from \eqref{3.3e} and replacing $q^2$ by $q$, we obtain
		\begin{equation}\label{3.3f}
			\sum_{n=0}^{\infty}\overline a_{2m+2}(4n+1)q^n\equiv2\dfrac{f_2^3f_4^5}{f_1^4f_{8}^2}\pmod{8}.
		\end{equation}
		Employing \eqref{yp2} in \eqref{3.3f}, we obtain
			\begin{equation}\label{3.3g}
				\sum_{n=0}^{\infty}\overline a_{2m+2}(4n+1)q^n\equiv2f_2f_4\pmod{8}.
			\end{equation}
    	Extracting the terms involving $q^{2n}$ from \eqref{3.3g} and replacing $q^2$ by $q$, we obtain
			\begin{equation}\label{3.3h}
				\sum_{n=0}^{\infty}\overline a_{2m+2}(8n+1)q^n\equiv2f_1f_2\pmod{8},
			\end{equation}
			which is the $\alpha=0$ case of \eqref{3.3}. Suppose that \eqref{3.3} is true for $\alpha\geq0$. Employing \eqref{u10} in \eqref{3.3}, we obtain
			$$\hspace{-10cm}\sum_{n=0}^{\infty}\overline a_{2m+2}(8\cdot p^{2\alpha}n+p^{2\alpha})q^n$$
			$$\equiv2\bigg[\sum_{\substack{ k={-(p-1)/2} \\ k \ne {(\pm p-1)/6}}}^{(p-1)/2}(-1)^{k}{q}^{(3k^2+k)/2} F\left( {-q}^{(3p^2+(6k+1)p)/2},{-q}^{(3p^2-(6k+1)p)/2}\right)+(-1)^{{(\pm p-1)/6}}{q}^{(p^2-1)/24}f_{p^2}\bigg]$$
			$$\hspace{.5cm}\times\bigg[\sum_{\substack{ k={-(p-1)/2} \\ k \ne {(\pm p-1)/6}}}^{(p-1)/2}(-1)^{k}{q}^{3m^2+m} F\left( {-q}^{3p^2+(6m+1)p},{-q}^{3p^2-(6m+1)p}\right)+(-1)^{{(\pm p-1)/6}}{q}^{(p^2-1)/12}f_{2p^2}\bigg] $$
			\begin{equation}\label{3.3j}
			\hspace{-12cm}\pmod{8}.
			\end{equation}
			Consider the congruence
			$$\dfrac{(3k^2+k)}{2}+\dfrac{2(3m^2+m)}{2}\equiv\dfrac{3(p^2-1)}{24}\pmod{p},$$
			which is equivalent to
			$$(6k+1)^2+2(6m+1)^2\equiv0\pmod{p}.$$
			Since $\left(\dfrac{-2}{p}\right)=-1$, the above congruence has only solution $k=m=\dfrac{(\pm p-1)}{6}$.\\
			Therefore, extracting the terms involving $q^{pn+(p^2-1)/8}$ from both sides of \eqref{3.3j}, dividing throughout by $q^{(p^2-1)/8}$ and then replacing $q^p$ by $q$, we obtain
			\begin{equation}\label{3.3l}
				\sum_{n=0}^{\infty}\overline a_{2m+2}(8\cdot p^{2\alpha+1}n+p^{2\alpha+2})q^n\equiv2f_pf_{2p}\pmod{8}.
			\end{equation}
			Extracting the terms involving $q^{pn}$ from \eqref{3.3l} and replacing $q^p$ by $q$, we obtain
			\begin{equation}\label{3.3m}
				\sum_{n=0}^{\infty}\overline a_{2m+2}(8\cdot p^{2(\alpha+1)}n+p^{2(\alpha+1)})q^n\equiv2f_1f_{2}\pmod{8},
			\end{equation}
				which is $\alpha+1$ case of \eqref{3.3}. Thus, by the principle of mathematical induction, the proof of \eqref{3.3} complete.\\
				Comparing the coefficients of the terms involving $q^{pn+r}$, for $1\leq r\leq p-1$, from both sides \eqref{3.3l}, we arrive at desired result \eqref{3.3i}.\end{proof}
	
	\begin{theorem}
		Let  $p\geq5$ be a prime with $\left(\dfrac{-8}{p}\right)=-1$ and $1\leq r \leq p-1$. Then for all integer $m\geq0$ and $\alpha\geq0$, we have
		\begin{equation}\label{3.4}
			\sum_{n=0}^{\infty}\overline a_{2m+2}(8\cdot p^{2\alpha}n+3\cdot p^{2\alpha})q^n\equiv4f_1f_8\pmod{8},
		\end{equation}and
		\begin{equation}\label{3.4i}
			\overline a_{2m+2}(8\cdot p^{2\alpha+1}(pn+r)+3\cdot p^{2\alpha+2})\equiv0\pmod{8}.
	\end{equation}\end{theorem}
	 \begin{proof}
	 	Employing \eqref{phim} in \eqref{3.3d}, we obtain
	 	\begin{equation}\label{3.4b}
	 		\sum_{n=0}^{\infty}\overline a_{2m+2}(2n+1)q^n\equiv2\dfrac{f_4^3}{\phi(-q)}\pmod{8}.
	 	\end{equation}
	 	Employing \eqref{1/phi} in \eqref{3.4b}, we obtain
	 	\begin{equation}\label{3.4c}
	 		\sum_{n=0}^{\infty}\overline a_{2m+2}(2n+1)q^n\equiv2\dfrac{f_4^3}{\phi(-q^4)^4}\left(\phi(q^4)^3+2q\phi(q^4)^2\psi(q^8)\right)\pmod{8}.
	 	\end{equation}
	 	Extracting the terms involving $q^{4n+1}$ from \eqref{3.4c}, dividing by $q$ and replacing $q^4$ by $q$ we obtain
	 	\begin{equation}\label{3.4d}
	 		\sum_{n=0}^{\infty}\overline a_{2m+2}(8n+3)q^n\equiv4\dfrac{f_1^3}{\phi(-q)^4}\left(\phi(q)^2\psi(q^2)\right)\pmod{8}.
	 	\end{equation}
	 	Employing \eqref{e0}, \eqref{tp1}, \eqref{phim} and \eqref{yp} in \eqref{3.4d}, we obtain
	 	\begin{equation}\label{3.4e}
	 			\sum_{n=0}^{\infty}\overline a_{2m+2}(8n+3)q^n\equiv4f_1f_8\pmod{8},
	 	\end{equation}
	 	which is the $\alpha=0$ case of \eqref{3.4}. Suppose that \eqref{3.4} is true for $\alpha\geq0$. Employing \eqref{u10} in \eqref{3.4}, we obtain
	 	$$\hspace{-10cm}\sum_{n=0}^{\infty}\overline a_{2m+2}(8\cdot p^{2\alpha}n+3\cdot p^{2\alpha})q^n$$
	 	$$\equiv4\bigg[\sum_{\substack{ k={-(p-1)/2} \\ k \ne {(\pm p-1)/6}}}^{(p-1)/2}(-1)^{k}{q}^{(3k^2+k)/2} F\left( {-q}^{(3p^2+(6k+1)p)/2},{-q}^{(3p^2-(6k+1)p)/2}\right)+(-1)^{{(\pm p-1)/6}}{q}^{(p^2-1)/24}f_{p^2}\bigg]$$
	 	$$\hspace{.5cm}\times\bigg[\sum_{\substack{ k={-(p-1)/2} \\ k \ne {(\pm p-1)/6}}}^{(p-1)/2}(-1)^{k}{q}^{4(3m^2+m)} F\left( {-q}^{4(3p^2+(6m+1)p)},{-q}^{4(3p^2-(6m+1)p)}\right)+(-1)^{{(\pm p-1)/6}}{q}^{(p^2-1)/3}f_{8p^2}\bigg] $$
	 	\begin{equation}\label{3.4f}
	 		\hspace{-12cm}\pmod{8}.
	 	\end{equation}
	 	Consider the congruence
	 	$$\dfrac{(3k^2+k)}{2}+4(3m^2+m)\equiv\dfrac{3(p^2-1)}{8}\pmod{p},$$
	 	which is equivalent to
	 	$$(6k+1)^2+8(6m+1)^2\equiv0\pmod{p}.$$
	 	Since $\left(\dfrac{-8}{p}\right)=-1$, the above congruence has only solution $k=m=\dfrac{(\pm p-1)}{6}$.\\
	 	Therefore, extracting the terms involving $q^{pn+3(p^2-1)/8}$ from both sides of \eqref{3.3j}, dividing throughout by $q^{3(p^2-1)/8}$ and then replacing $q^p$ by $q$, we obtain
	 	\begin{equation}\label{3.4g}
	 		\sum_{n=0}^{\infty}\overline a_{2m+2}(8\cdot p^{2\alpha+1}n+3\cdot p^{2\alpha+2})q^n\equiv4f_pf_{8p}\pmod{8}.
	 	\end{equation}
	 	Extracting the terms involving $q^{pn}$ from \eqref{3.4g} and replacing $q^p$ by $q$, we obtain
	 	\begin{equation}\label{3.4h}
	 		\sum_{n=0}^{\infty}\overline a_{2m+2}(8\cdot p^{2(\alpha+1)}n+3\cdot p^{2(\alpha+1)})q^n\equiv4f_1f_{8}\pmod{8},
	 	\end{equation}
	 	which is $\alpha+1$ case of \eqref{3.4}. Thus, by the principle of mathematical induction, the proof of \eqref{3.4} complete.\\
	 	Comparing the coefficients of the terms involving $q^{pn+r}$, for $1\leq r\leq p-1$, from both sides \eqref{3.3l}, we arrive at desired result \eqref{3.4i}.
	 \end{proof}
	\begin{corollary}For any integers $m$, $n\geq0$, we have
		\begin{equation}\label{c6}
			a_{2m+2}(8n+z)\equiv0\pmod{8},\quad for \quad z=5,7.
		\end{equation}
		\begin{proof}
				Extracting the terms involving $q^{4n+i}$ for $i=2,3$ from \eqref{3.4c}, we easily arrived at \eqref{c6}.
		\end{proof}
	\end{corollary}
	
	\begin{theorem}
		Let $p\geq3$ is any prime. Then for any integer  $\alpha\geq0$, $m\geq0$, $n\geq0$ and $1\leq r\leq p-1$, we have
		\begin{equation}\label{3.5}
			\sum_{n=0}^{\infty}\overline a_{8m+3}(8\cdot p^{2\alpha}n+p^{2\alpha})q^n\equiv2\psi(q)\pmod{16},
		\end{equation}and
		\begin{equation}\label{3.5i}
			\overline a_{8m+3}(8\cdot p^{2\alpha+1}(pn+r)+p^{2\alpha+2})\equiv0\pmod{16}.
	\end{equation}\end{theorem}
	
	\begin{proof}
			Setting $c=3$ in \eqref{goc}, we obtain
		\begin{equation}\label{3.5a}
			\sum_{n=0}^{\infty}\overline a_3(n)q^n=\dfrac{f_4^2}{f_1^2f_2^3}.
		\end{equation}
		Employing \eqref{f_1^2} in \eqref{3.5a}, we obtain
		\begin{equation}\label{3.5b}
			\sum_{n=0}^{\infty}\overline a_3(n)q^n=\dfrac{f_4^2f_8^5}{f_2^8f_{16}^2}+2q\dfrac{f_4^4f_{16}^2}{f_2^8f_8}.
		\end{equation}
		Extracting the terms involving $q^{2n+1}$ from \eqref{3.5b}, dividing by $q$ and replacing $q^2$ by $q$ we obtain
		\begin{equation}\label{3.5c}
			\sum_{n=0}^{\infty}\overline a_3(2n+1)q^n=2\dfrac{f_2^4f_{8}^2}{f_1^8f_4}.
		\end{equation}
		Employing \eqref{yp2} in \eqref{3.5c} and then applying \eqref{3.1}, we obtain
		\begin{equation}\label{3.5d}
			\sum_{n=0}^{\infty}\overline a_{8m+3}(2n+1)q^n\equiv2\dfrac{f_8^2}{f_4}\pmod{16}.
		\end{equation}
		Extracting the terms involving $q^{4n}$ from \eqref{3.5d} and replacing $q^4$ by $q$, we obtain
		\begin{equation}\label{3.5e}
				\sum_{n=0}^{\infty}\overline a_{8m+3}(8n+1)q^n\equiv2\dfrac{f_2^2}{f_1}\pmod{16}.
		\end{equation}
		 Employing \eqref{tp1} in \eqref{3.5e}, we obtain
		   \begin{equation}\label{3.5f}
		   		\sum_{n=0}^{\infty}\overline a_{8m+3}(8n+1)q^n\equiv2\psi(q)\pmod{16},
		   \end{equation} 
		 which is the $\alpha=0$ case of \eqref{3.5}. Suppose that \eqref{3.5} is true for $\alpha\geq0$. Employing \eqref{w4} in \eqref{3.5}, we obtain
		    $$\hspace{-10cm}\sum_{n=0}^{\infty}\overline a_{8m+3}(8\cdot p^{2\alpha}n+ p^{2\alpha})q^n$$
		    \begin{equation}\label{3.5g}
		    	\equiv2\sum_{j=0}^{(p-3)/2}{\it q}^{(j^2+j)/2} F\left( {\it q}^{(p^2+(2j+1)p)/2},{\it q}^{(p^2-(2j+1)p)/2}\right)+{\it q}^{(p^2-1)/8}\psi(q^{p^2})\pmod{16}.
		    \end{equation}
		    Extracting the terms involving $q^{pn+(p^2-1)/8}$ from both sides of \eqref{3.5g}, dividing throughout by $q^{(p^2-1)/8}$ and then replacing $q^p$ by $q$, we obtain
		    	\begin{equation}\label{3.5h}
		    	\sum_{n=0}^{\infty}\overline a_{8m+3}(8p^{2\alpha+1}n+p^{2\alpha+2})q^n\equiv2\psi(q^p)\pmod{16}.
		    \end{equation}
		    Extracting the terms involving $q^{pn}$ from \eqref{3.5h} and replacing $q^p$ by $q$, we obtain
		    \begin{equation}\label{3.5j}
		    	\sum_{n=0}^{\infty}\overline a_{8m+3}(8p^{2(\alpha+1)}n+p^{2(\alpha+1)})q^n\equiv2\psi(q)\pmod{8},
		    \end{equation}
		    which is $\alpha+1$ case of \eqref{3.5}. Thus, by the principle of mathematical induction, the proof of \eqref{3.5} complete.\\
		    Comparing the coefficients of the terms involving $q^{pn+r}$, for $1\leq r\leq p-1$, from both sides \eqref{3.5h}, we arrive at desired result \eqref{3.5i}.	\end{proof}
		    
		    \begin{corollary}For any integers $m$, $n\geq0$, we have
		    	\begin{equation}\label{c5}
		    		a_{8m+3}(8n+w)\equiv0\pmod{16}, \quad for \quad w=3,5,7.
		    	\end{equation}
		    	\begin{proof}
		    			Extracting the terms involving $q^{4n+i}$ for $i=1,2,3$ from \eqref{3.5d}, we easily arrived at \eqref{c5}.
		    	\end{proof}
		    \end{corollary}

\end{document}